\documentclass[10pt,reqno,a4paper]{amsart}

\usepackage[pdftex]{graphicx}
\usepackage{a4wide}
\usepackage{amsmath}
\usepackage{amssymb,amsthm,hyperref,caption}
\usepackage{color}
\definecolor{blau}{rgb}{0,0,0.75} 
\hypersetup{colorlinks,linkcolor=blau,citecolor=blue}

\allowdisplaybreaks

\newtheorem{theorem}{Theorem}
\newtheorem{lemma}{Lemma}
\newtheorem{coroll}{Corollary}

\theoremstyle{definition}
\newtheorem{remark}{Remark}

\def\P{{\mathbb {P}}}

\newcommand{\N}{\ensuremath{\mathbb{N}}}

\newcommand{\R}{\ensuremath{\mathbb{R}}}
\newcommand{\Z}{\ensuremath{\mathbb{Z}}}

\begin{document}

\author{Markus Kuba}
\address{Markus Kuba\\
Institut f{\"u}r Diskrete Mathematik und Geometrie\\
Technische Universit\"at Wien\\
Wiedner Hauptstr. 8-10/104\\
1040 Wien, Austria} \email{kuba@dmg.tuwien.ac.at}


\title{A note on naturally embedded ternary trees}

\begin{abstract}
In this note we consider ternary trees naturally embedded in the plane in a deterministic way such that the root has position zero, or in other words label zero, and the children of a node with position $j$ have positions $j-1$, $j$, and $j+1$, for all $j\in\Z$. 
We derive the generating function of ternary trees where all nodes have labels which are less or equal than $j$, with $j\in\N$, which 
generalizes a result of~\cite{Schaeff1998} and~\cite{Left2000}, and the generating function of ternary trees counted with respect to nodes with label $j$, with $j\in\Z$. Moreover, we discuss generalizations of the counting problem to several labels at the same time. Furthermore, we use generating functions to study the depths of the external node $s$, or in other words leaf $s$ with $0\le s\le 2n$, where the $2n+1$ external nodes of a ternary tree are numbered from the left to the right according to an inorder traveral. The three different types depths -- left, right and center -- are due to the embedding of the ternary tree in the plane. 
Finally, we discuss generalizations of the considered enumeration problems to embedded $d$-ary trees. 
\end{abstract}

\thanks{The author was supported by the Austrian Science Foundation FWF, grant S9608-N13.}
\keywords{Ternary trees, Embedded trees, Labelled trees}%

\maketitle
\section{Introduction}
The study of tree families embedded in the plane has recently received a lot of attention. Binary trees, complete binary trees, more generally simply generated tree families, and several different families of planar trees have been considered in a series of papers~\cite{Boutt2003,Boutt2003II,Marck2004,Francesco2005,Bou2006,BouJan2006,Sva2006,Kne2006,DevSva2008,Alois2009}. 
It has been shown that embedded trees are closely related to a random measure called
ISE (Integrated SuperBrownian Excursion). For example,
in the recent paper of Devroye and Janson~\cite{DevSva2008} a
conjecture of Bousquet-M\'elou and Janson~\cite{BouJan2006} is
proven, saying that the vertical profile of a randomly labelled
simply generated tree converges in distribution, after
suitable normalization, to the density of the ISE. However, enumerative properties of deterministically embedded trees 
have not been intensively studied, except for binary trees, a particular subclass of ternary trees, and families of plane trees. Motivated by the results
of Bousquet-M\'elou~\cite{Bou2006} and Panholzer~\cite{Alois2009}
for embedded binary trees, and the results of Schaeffer and Jacquard~\cite{Schaeff1998}, and Del Lungo, Del Ristoro and Penaud~\cite{Left2000} for a specific subclass of embedded ternary trees, we study several combinatorial properties
of ternary trees embedded in the plane in a \emph{deterministic}
manner. We consider the following natural embedding. The root has
position zero, or in other words label zero, and
the three children of the root have positions $-1$, $0$ and $1$.
More generally, the labels of the children of a node with
label $j$ are given by $j-1,j,j+1$, with $j\in\Z$. 
A similar embedding of ternary trees has been considered before~\cite{Schaeff1998,Left2000}, 
where the authors studied a particular subclass of embedded ternary trees named skew ternary trees~\cite{Schaeff1998}, or left ternary trees~\cite{Left2000}, which are embedded ternary trees with no node having label greater than zero. Using bijections between embedded ternary trees with no label greater than zero and non-separable rooted planar maps with $n + 1$ edges they obtained amongst others an explicit result for the number of such trees of size $n$. The aim
of this work is to use a generating functions approach to study
several parameters in embedded ternary trees: we analyze the
distribution of the different types of depths of external nodes in
ternary trees, stemming from the embedding of the tree, assuming
that the external nodes are enumerated from the left to the right, or in other words with respect to an inorder-traversal. 
Related parameters have been studied for binary trees
in~\cite{PP1997,PP2002,Alois2009}. Moreover, we are interested in
the number of embedded ternary trees of size $n$ where all internal
nodes have label smaller or equal than $j$, with $j\in\N$, where we extend the results of~\cite{Schaeff1998} and~\cite{Left2000} for the special case $j=0$ to arbitrary $j\ge 0$,  
and also in the number of embedded ternary of size $n$ counted with
respect to the number of internal nodes with label $j$,
with $j\in\Z$. We also show how the extend the counting problem to
several types of labels considering the nodes with label $j$ and the number of nodes with labels in $\{j-1,j+1\}$, and
also discuss generalizations. 
\subsection{Plan of the paper} This note is structured as
follows. In the next section we recall some well known properties of the family
of ternary trees. In Section~\ref{DEEBsecembedd} we discuss the natural embedding of
(ternary) trees in the plane, and also provide a formalism for embedded trees.
Section~\ref{DEEBboutt} is devoted to the study of ternary trees with small label and to ternary
trees counted by the number of nodes labelled $j$. 
In Subsection~\ref{DEEBSubsectiongen1} we discuss the enumeration
problem of counting the number of nodes with label $j$, and at the
same time the number of nodes with labels in $\{j-1,j+1\}$. In
Section~\ref{DEEBTERleaves} we study the distribution of depths of
the $s$-th external node or leaf in ternary trees, with $0\le s\le 2n$, where the
leaves are enumerated from the left to the right. In the final
section we discuss some generalizations and open problems, i.e.~the extensions of the
obtained results to embedded $d$-ary trees.


\section{Preliminaries}
\subsection{The family of ternary trees}
The family of ternary trees $\mathcal{T}$ can be described in a recursive way, which
says that a ternary tree is either a leaf (an external node) or an internal node followed by three ordered ternary trees, visually described by
the suggestive ``equation''
\begin{align*}
\label{DEEBaloissymb1}
\parbox{4cm}{
\unitlength 0.50mm
\linethickness{0.4pt}
\thicklines
\begin{picture}(78.33,22.00)
\put(64.33,1.67){\circle{6.57}}
\put(64.33,4.67){\line(0,1){10}}
\put(62.00,4.00){\line(-1,1){11}}
\put(66.67,4.00){\line(1,1){11}}
\put(42.00,14.33){\makebox(0,0)[cc]{+}}
\put(21.00,13.50){\makebox(0,0)[cc]{\large{$\Box$}}}
\put(50.33,20.00){\makebox(0,0)[cc]{$\mathcal{T}$}}
\put(64.33,20.00){\makebox(0,0)[cc]{$\mathcal{T}$}}
\put(78.33,20.00){\makebox(0,0)[cc]{$\mathcal{T}$}}
\put(1.67,14.67){\makebox(0,0)[cc]{$\mathcal{T}\quad=\quad$}}
\end{picture}}
\end{align*}
Here $\text{\small{$\bigcirc$}}$ is the symbol for an internal node and
$\text{\large{$\Box$}}$ is the symbol for a leaf or external node.
A trivial consequence of this description is the fact that a ternary tree with $n$ internal nodes has exactly $2n+1$ external nodes;
moreover, by taking external nodes into account any node has either outdegree zero or three. Note that throughout this work the notion size of a tree is with respect to the number of internal nodes, a tree of size $n$ has $n$ internal nodes. We assume that the $2n+1$ external nodes of a size $n$ ternary tree are numbered from left to right according to a so-called inorder traversal. We start at the root node of a tree. If the tree has internal nodes, we recursively traverse them by going first to the left subtree, then the center subtree, and finally to the right subtree. The external nodes are numbered as visited on the traversal process, see Figure~\ref{DEEBfig2}.
The generating function $T(z)=\sum_{n\ge 0}T_n z^n$ of the number of ternary trees of size $n$ satisfies the equation
\begin{equation}
\label{DEEBteReqn0yi}
T(z)=1+zT^3(z),\quad\text{with}\,\,T(0)=1.
\end{equation}
Concerning the series expansion of the generating function $T(z)$ it is convenient consider the shifted series $\tilde{T}(z):=T(z)-1$. This corresponds to discarding external nodes (the empty tree) in the discription above; we obtain simply generated ternary trees $\tilde{\mathcal{T}}$ ), defined by the formal equation
\begin{equation}
   \label{DEEBteReqn1}
   \tilde{\mathcal{T}} = \bigcirc \times \Big(1 \cdot \{\epsilon\} \; \dot{\cup} \;
   3 \cdot \tilde{\mathcal{T}} \; \dot{\cup} \; 3 \cdot
   \tilde{\mathcal{T}} \times \tilde{\mathcal{T}}
   \; \dot{\cup} \; 1 \cdot \tilde{\mathcal{T}} \times \tilde{\mathcal{T}} \times \tilde{\mathcal{T}}\;\Big)
   = \bigcirc \times \varphi(\tilde{\mathcal{T}}),\quad\text{with}\,\, \varphi(t)=(1+t)^3,
\end{equation}
with $\bigcirc$ a node, $\times$ the cartesian product, and
$\varphi(\tilde{\mathcal{T}})$ the substituted structure. We refer to~\cite{MeirMoon1978} for the general definition of simply generated trees.
Let $T_n$ denote the number of ternary trees of size $n$, and $\tilde{T}_n$ the number of simply generated ternary trees of size $n$. By the formal description above~\eqref{DEEBteReqn1} the counting series $\tilde{T}(z)=\sum_{n\ge 1}\tilde{T}_nz^n$ satisfies
the functional equation
\begin{equation}
\label{DEEBteReqn2}
\tilde{T}(z)=z(1+\tilde{T}(z))^3,\quad \tilde{T}(0)=0.
\end{equation}
Due to the Lagrange inversion formula, see e.g.~\cite{Gould1983}, the number of ternary trees of size $n$ is given by
\begin{equation}
\label{DEEBteReqn3}
\tilde{T}_n=[z^n]\tilde{T}(z)=\frac1{2n+1}\binom{3n}{n},\quad\text{and consequently}\quad \tilde{T}(z)=\sum_{n\ge 1}\binom{3n}{n}\frac{z^n}{2n+1}.
\end{equation}
Note that due to the definition the series $T(z)$ and $\tilde{T}(z)$ are related by $T(z)=\tilde{T}(z)-1$.
In contrast to the case of binary trees, where the generating function $y(z)=\sum_{n\ge 1}\binom{2n}{n}\frac{z^n}{n+1}$, satisfying $y=z(1+y)^2$, may be expressed in terms of radicals, $y=(1-2z-\sqrt{1-4z})/(2z)$, for $|z|\le 1/4$, the case of ternary trees is more complicated due to the third order equation~\eqref{DEEBteReqn2}; however, there exist appealing closed form expressions of $\tilde{T}(z)$ for real values of $z$, i.e.~$z\in[-\frac4{27},0]$ and $z\in[0,\frac4{27}]$, 
\begin{equation*}
\tilde{T}(z)=
\begin{cases}
\frac{2}{\sqrt{3z}}\sin\Big(\frac13\arctan\big(\sqrt{\frac{27z}{4-27z}}\big)\Big)-1, &\,\,\text{for}\,\,z\in[0,\frac4{27}],\\[0.3cm]
\frac{\sqrt[3]{4(\sqrt{4-27z}+\sqrt{-27z})}-\sqrt[3]{4(\sqrt{4-27z}-\sqrt{-27z})}}{2\sqrt{-3z}}-1,&\,\,\text{for}\,\,z\in[-\frac4{27},0],
\end{cases}
\end{equation*}
which are obtained by Cardano's method; note that the radius of convergence of $\tilde{T}(z)$ and $T(z)$ is given by $4/27$. We obtain by the Lagrange inversion formula the more general formula
\begin{equation*}
[z^n]\big(\tilde{T}(z)\big)^k=\frac{k}{n}\binom{3n}{n-k},\quad\text{for}\,\,k\in\N.
\end{equation*}
Note that since $T(z)=\tilde{T}(z)-1$, this immediately implies that for $n\ge k$ we have
\begin{equation}
\label{DEEBteReqn4}
[z^n]\big(T(z)\big)^{k}=[z^n]\big(\tilde{T}(z)+1\big)^k=\sum_{\ell=1}^{k}\binom{k}{\ell}\frac{\ell}{n}\binom{3n}{n-\ell}
=\frac{k}{n}\binom{3n+k-1}{n-1}=\frac{k}{3n+k}\binom{3n+k}{n}.
\end{equation}
It turns out that the last representation is valued for all $n\ge
0$. The result above follows by the binomial theorem,
Equation~\ref{DEEBteReqn3} and a variant of the Chu-Vandermonde
summation formula.

\section{The embedding of ternary trees\label{DEEBsecembedd}}
\begin{figure}[htb!]
\begin{minipage}[b]{0.35\linewidth}
\centering
\includegraphics[angle=0,scale=0.8]{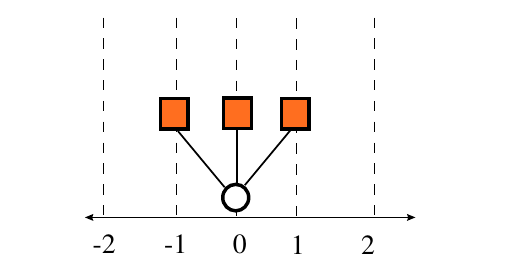}
\caption{An embedded size 1 ternary tree with its external nodes.}
\label{DEEBfig1}
\end{minipage}
\setcaptionwidth{0.55\textwidth}
\begin{minipage}[b]{0.6\linewidth}
\centering
\includegraphics[angle=0,scale=0.9]{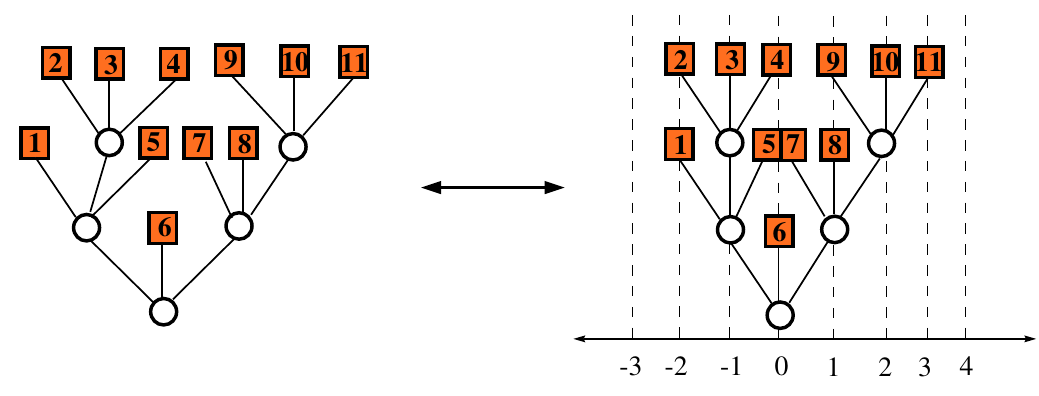}
\caption{An embedded ternary tree of size 5, its external nodes numbered w.~r.~t~the inorder-traversal.}
\label{DEEBfig2}
\end{minipage}
\end{figure}

Motivated by the results of Bousquet-M\'elou~\cite{Bou2006}, and Panholzer~\cite{Alois2009} for naturally embedded binary trees, 
and the results of Schaeffer and Jacquard~\cite{Schaeff1998}, and Del Lungo, Del Ristoro and Penaud~\cite{Left2000} for a subclass 
of ternary trees, we embedd ternary trees in the plane in the following way. The root node has position zero. By definition of ternary trees each internal node with no children has exactly three positions to attach a new node, which are as usual called external nodes or leaves, see Figure~\ref{DEEBfig1}. We recursively define the embedding of ternary increasing trees as follows: an internal node with label/position $j\in\Z$ has exactly three children, being internal or external, placed at positions $j-1$, $j$ and $j+1$. Following~\cite{Bou2006}, we call this embedding the \emph{natural embedding} of ternary trees, because the label a node is its abscissa in the natural integer embedding of the tree. 
\setcaptionwidth{0.3\textwidth}

\begin{figure}[htb]
\centering
\includegraphics[angle=0,scale=0.8]{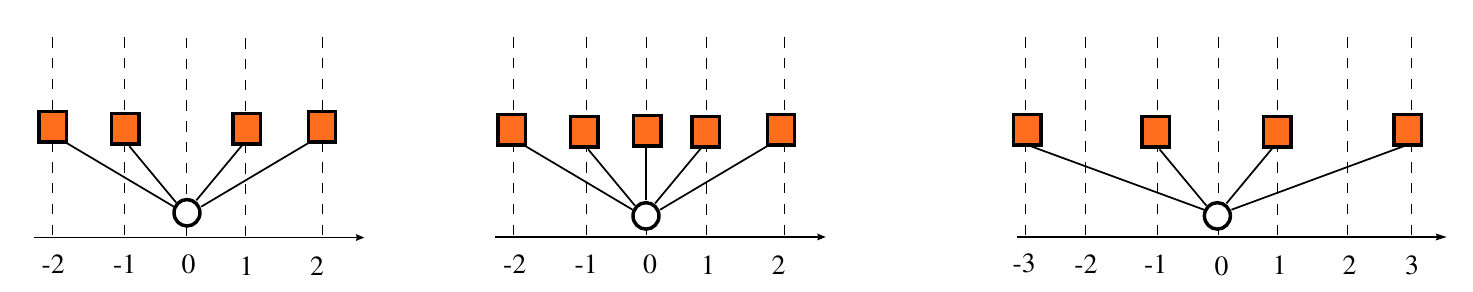}
\caption{Size one embedded quaternary and quinary trees together with their external nodes, i.e.~the possible increments $i\in\{\pm 1, \pm 2\}$ and $i\in\{0\}\cup\{\pm 1,\pm 2\}$; an alternative model for embedded quinary trees with increments $i\in\{\pm 1, \pm 3\}$.}
\label{DEEBfig3}
\end{figure}

\vspace{-0.6cm}
Equivalent to the embedding, we can think of the edges $e$ of a ternary tree being decomposed into three different types, $e_1$, $e_2$ and $e_3$.
We assume that each internal nodes $v$ has exactly one of each edge type $e_i$ pointing away from $v$.
Hence, if we assume that an edge of type $i$, with $1\le i\le 3$, corresponds to a step $2-i$, the label $l(v)$
of a node $v$ in a tree $B$ corresponds to the sum over all edges on the unique path $P=P(B)$ from node $v$ to the root,
\begin{equation*}
l(v)=\sum_{e_1\in P(B)}(-1) +\sum_{e_2\in P(B)}0 + \sum_{e_3\in P(B)}1= \sum_{e_3\in P(B)}1-\sum_{e_1\in P(B)}1.
\end{equation*}
Note that we associate two different kind of labels to external nodes, first the number of the external node according to the inorder-traversal, and second the label of its position according to the natural embedding. 
This definition readily extends to $d$-ary trees; following the embedding of Bousquet-M\'elou~\cite{Bou2006}
of binary trees and the present one stated for ternary trees. It seems natural to assume that for $(2d)$-ary trees
each internal node with label $j\in\Z$ has exactly $2d$ children, internal or external, placed
at positions $j+\ell$, with $\ell\in\{\pm 1,\dots, \pm d\}$, and for $(2d+1)$-ary trees
each internal node with label $j\in\Z$ has exactly $2d+1$ children, internal or external, placed
at positions $j+\ell$, with $\ell\in\{0\}\cup\{\pm 1,\dots, \pm d\}$; see Figure~\ref{DEEBfig3}.
Equivalently, we have $(2d)$ or $(2d+1)$ different types
of edges $e_{\ell}$, with either $\ell\in\{\pm 1,\dots, \pm d\}$ or $\ell\in\{0\}\cup\{\pm 1,\dots, \pm d\}$ which correspond to the steps in the obvious way. We point out that the terminology ``natural embedding'' may be not unique. One could alternatively define the embedding of $(2d)$-ary trees by the increments $\ell\in\{\pm 1,\pm 3\dots, \pm (1+2d)\}$, which could be considered more natural than $\ell\in\{\pm 1,\dots, \pm d\}$, depending on the taste.
\setcaptionwidth{0.8\textwidth}

\vspace{-0.6cm}

\subsection{A formalism for embedded trees\label{DEEBform}}
Previously we have stated a definition for embedded ternary trees and more general models of embedded $(2d+1)$-~and $(2d)$-ary trees. However, we have not presented a formalism to actually encode a specific embedded (ternary) tree. Moreover, there exist several different ``natural'' embeddings of $d$-ary trees. In the following we will give a simple natural way of encoding families of embedded trees, with respect to a set of step vectors. Let $\mathcal{S}$ denote a given finite non-empty set of step vectors $\mathcal{S}=\{\mathbf{a}_1,\dots,\mathbf{a}_d\}$ of size $|\mathcal{S}|=d$, with $\mathbf{a}_{\ell}=(1,b_{\ell})$ and $b_{\ell}\in\Z$, for $1\le \ell \le d$; moreover let $A$ denote the alphabet $A=\{1,\dots,d\}$. An \emph{embedded tree} $\mathcal{M}=\mathcal{M}(\mathcal{S})$ of size one, $\mathcal{M}=\{\mathbf{m}_1\}=\{(0,0,\epsilon)\}$, is defined as a single point located as the origin with an additional label $\epsilon$, denoting the empty word, where the first two entries of the triple $\mathbf{m}_1$ encode the points coordinate, and the third entry is the label.
An embedded tree $\mathcal{M}=\mathcal{M}(\mathcal{S})$ of size $n\ge 1$ is defined as a set of labelled points $\mathcal{M}=\{\mathbf{m}_1,\dots,\mathbf{m}_n\}$, satisfying two conditions stated below, with $\mathbf{m}_{\ell}$ denoting a triple $\mathbf{m}_{\ell}=( \mathbf{x}_{\ell},\lambda_{\ell})$, where $\mathbf{x}_{\ell}=(x_{\ell},y_{\ell})$ denote the coordinates of point $\mathbf{m}_{\ell}$ and $\lambda_{\ell}\in A^{*}$ denote its label consisting of a word of length $r$ over the alphabet $A$, with $0\le r\le n-1$. We impose two conditions on embedded trees. First, we impose that $(0,0,\epsilon)\in \mathcal{M}$, and that no other point $\mathbf{m}_{\ell}=( \mathbf{x}_{\ell},\lambda_{\ell})\in\mathcal{M}$ has a label of length zero. Second, for each labelled point $\mathbf{m}_{\ell}= (\mathbf{x}_{\ell},\lambda_{\ell})$, with $\lambda\neq \epsilon$, there exists an index $\ell'\in\{1,2\dots,n\}$ with corresponding point $\mathbf{m}_{\ell'}=(\mathbf{x}_{\ell'},\lambda_{\ell'})$ and an index $i$, with $1\le i \le d$, such that the point $\mathbf{m}_{\ell}=(\mathbf{x}_{\ell},\lambda_{\ell})$ is given by $\mathbf{m}_{\ell} =(\mathbf{x}_{\ell'}+\mathbf{a}_{i},\lambda_{\ell'}i)$.
Equivalently, given an embedded tree $\mathcal{M}=\mathcal{M}(\mathcal{S})=\{\mathbf{m}_1,\dots,\mathbf{m}_n\}$ of size $n$, we may obtain a size $n+1$ embedded tree $\mathcal{M}'$ by adding the points of the form $\mathbf{m}_{n+1}=(\mathbf{x}_{\ell}+\mathbf{a}_i,\lambda_{\ell}i)$ to $\mathcal{M}$, with $1\le \ell \le n$ and $1\le i\le d$. 

\setcaptionwidth{0.9\textwidth}
\begin{figure}[htb]
\centering
\includegraphics[angle=0,scale=0.8]{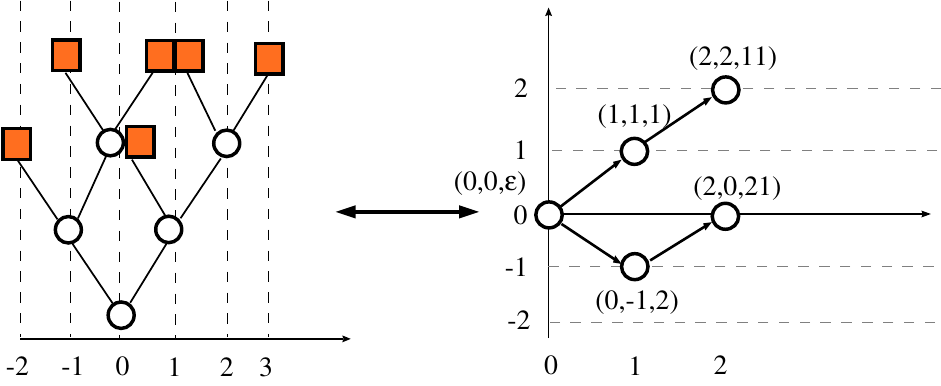}
\caption{An embedded binary tree of size $5$ and its corresponding size $5$ path configuration $\mathcal{M}=\{(0,0,\epsilon),(1,1,1),(1,-1,2),(2,0,12),(2,-2,22)\}$, with $\mathcal{S}=\{\mathbf{a}_1,\mathbf{a}_2\}$, $\mathbf{a}_1=(1,1)$ and $\mathbf{a}_2=(1,-1)$.\label{Deebpath}}
\end{figure}

\smallskip

Informally speaking, the stated description of embedded trees corresponds to a reflection followed by a 90 degree rotation of the usual pictures of embedded tree, see Figure~\ref{Deebpath}. Since there may be multiple points at the same coordinate $\mathbf{x}$, we use their labels to distinguish between them; more precisely the labels of the points encode how they are reached from the origin with respect to the given set of step vectors $\mathcal{S}=\{\mathbf{a}_1,\dots,\mathbf{a}_d\}$. The first condition says nothing else that all trees are rooted at the origin, and the second imposes that each non-root node has a unique parent node. The set of vectors $\mathcal{S}=\{(1,1),(1,0),(1,-1)\}$, i.e.~the set of vectors of Motzkin paths~\cite{BandFla2002}, leads to embedded ternary trees; $\mathcal{S}=\{(1,1),(1,-1)\}$, i.e.~the set of vectors of Dyck paths~\cite{BandFla2002}, leads to embedded binary trees as considered in~\cite{Bou2006}.
More generally, the family of embedded $(2d+1)$-ary and $(2d)$-ary trees, with respect to the natural embedding such that each internal node with label $j$ has exactly $(2d+1)$ and $(2d)$ children, internal or external, placed at positions $j+\ell$, with $\ell\in\{0\}\cup\{\pm 1,\dots, \pm d\}$ and $\ell\in\{\pm 1,\dots, \pm d\}$, 
is described using the stepsets $\mathcal{S}_{2d+1}=\{(1,\ell) \mid \ell\in \{0\}\cup\{\pm 1,\dots, \pm d\}\}$, and $\mathcal{S}_{2d}=\{(1,\ell) \mid \ell \in \{\pm 1,\dots, \pm d\}\}$.

\smallskip

The family of naturally embedded ternary trees with $\mathcal{S}=\{(1,1),(1,0),(1,-1)\}$ and the family of naturally embedded binary trees with $\mathcal{S}=\{(1,1),(1,-1)\}$, considered in~\cite{Bou2006}, are \emph{symmetric}: if $(1,b_{\ell})\in\mathcal{S}$ then also $(1,-b_{\ell})\in\mathcal{S}$.

\smallskip

For a given set of step vectors $\mathcal{S}$ let $\mathcal{F}=\mathcal{F}(\mathcal{S})$ denote the corresponding family of embedded trees. It is easy to see that the generating function $T(z)=\sum_{\mathcal{M}\in\mathcal{F}}z^{|\mathcal{M}|}$ satisfies the equation $T(z)=1+zT^{|S|}(z)$, where we assume that there exists exactly one embedded tree of size zero. Consequently, with respect to enumeration, embedded trees with $|\mathcal{S}|=d$ can be considered as models of $d$-ary trees.


\section{The number of embedded ternary trees with small labels\label{DEEBboutt}}
Let $T_j(z)$ denote the generating function of ternary trees having no label greater than $j$, with $j\ge 0$,
and with $T(z)$ the generating function of ternary trees, as specified by~\eqref{DEEBteReqn0yi}.
The starting point of our considerations is the following system of equations.
\begin{lemma}
The series $T_j(z)$ satisfies
\begin{equation}
\label{BOUTT1}
T_j(z)=1+zT_{j-1}(z)T_j(z)T_{j+1}(z),\quad \text{for}\,\, j=0,1,2,\dots,\quad\text{with}\,\,T_{-1}(z)=1.
\end{equation}
\end{lemma}

\begin{proof}
First we observe that the infinite system of equations is well defined and complete determines the generating functions $T_j(z)$.
Moreover, following~\cite{Bou2006} we note that replacing each label $k$ by $j-k$ shows that
the series $T_j(z)$ is also the generating function of trees rooted at a node with label $j$,
and having only non-negative labels for its children. Considering such a tree
it has at most three subtrees rooted at $j+\ell$, with $\ell\in\{-1,0,1\}$, which have again only non-negative labels.
We have the formal description sketched below in Figure~\ref{DEEBfig4}, which translates in the stated system of equations
and the result follows. \end{proof}
\begin{figure}[!htb]
\centering
\includegraphics[angle=0,scale=0.8]{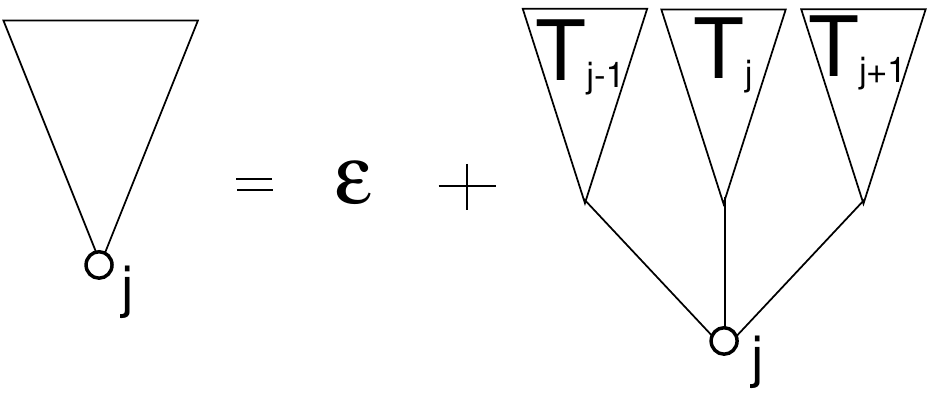}
\caption{The formal decomposition of embedded ternary trees.\label{DEEBfig4}}
\end{figure}

\vspace{-0.8cm}

We obtain the following result for the series $T_j(z)$.
\begin{theorem}
\label{DEEBtheBOUTT1} Let $T_j(z)$ be the generating function of
ternary trees with no label greater than $j$. Then $T_j(z)$ is given
by the the following expression
\begin{equation*}
T_j(z)=T(z)\frac{(1-X^{j+2}(z))(1-X^{j+5}(z))}{(1-X^{j+3}(z))(1-X^{j+4}(z))},\quad \text{for}\,\,j\ge -1,
\end{equation*}
where $T=T(z)$ is defined by~\eqref{DEEBteReqn0yi}, i.e.~$T=1+zT^3$, and the series
$X=X(z)$ is defined as the unique formal power series with $X(0)=0$, satisfying the relations
\begin{equation*}
X=zT^2X(\frac{1}{X}+1+X),\qquad X=z\frac{(1+X+X^2)^3}{(1+X^2)^2},\quad\text{and}\quad T=\frac{1+X+X^2}{1+X^2}.
\end{equation*}
Moreover, the series $X$ has non-negative coefficients and satisfies
\begin{equation*}
X=\frac{1-zT^2-\sqrt{1-2zT^2-3z^2T^4}}{2zT^2}=\frac{1-\sqrt{1-4\tilde{T}^2}}{2\tilde{T}},
\end{equation*}
with $\tilde{T}=\tilde{T}(z)=T(z)-1$, as defined in~\eqref{DEEBteReqn2}.
\end{theorem}
\begin{remark}
Let $U_n(w)$ denote the $n$-th Chebyshev polynomial of the second kind,
recursively defined by $U_0(z)=1$, $U_1(z)=2w$ and $U_{n+1}(w)=2wU_n(w)-U_{n-1}(w)$ for $n\ge 1$.
Following~Bouttier et al.~\cite{Boutt2003} we use the well known closed form expression of $U_n(w)$,
\begin{equation*}
U_n(w)=\frac{(w+\sqrt{w^2-1})^{n+1}-(w-\sqrt{w^2-1})^{n+1}}{2\sqrt{w^2-1}},
\end{equation*}
and observe that $T_j(z)$ can be expressed as the quotient of Chebyshev polynomials of the second kind evaluated at $w=(\sqrt{X}+1/\sqrt{X})/2$, we have for $j\ge 0$ the expression
\begin{equation*}
T_j(z)=T(z)\frac{U_{j+1}(w)U_{j+4}(w)}{U_{j+2}(w)U_{j+3}(w)}.
\end{equation*}
\end{remark}
An immediate consequence of our result above are explicit formulas for $T_j(z)$, i.e.~$T_0(z)$ and $T_1(z)$, stated in the following corollary, which is obtained by elimination of $X$ using standard algorithms from the theory of Gr\"obner bases.
\begin{coroll}
\label{DEEBcoroll}
The series $T_0(z)$ and $T_1(z)$, counting the number of embedded ternary trees with no label greater than zero, and one is given by following simple expression.
\begin{equation*}
T_0(z)=3T(z)-1-T^2(z),\quad T_1(z)=\frac{(T(z)-2)T^3(z)}{T^2(z)-3T(z)+1}.
\end{equation*}
The coefficients $[z^n]T_0(z)$ and $[z^n]T_1(z)$ are for $n\ge 1$ given by the expressions
\begin{equation*}
\begin{split}
[z^n]T_0(z)&=\frac{2}{(n+1)(2n+1)}\binom{3n}{n},\\
[z^n]T_1(z)&=\frac{2}{n+1}\binom{3n}{n}+\sum_{k=0}^{n}(-1)^{k+1}F_{k+1}\binom{3n}{n-k}
\frac{n(11k+5) -2k(k+1)}{n(2n+k+1)},
\end{split}
\end{equation*}
where $F_k=\frac{\phi^k- \frac{(-1)^k}{\phi^k}}{\sqrt{5}}$ denotes the $k$-th Fibonacci number, and $\phi=\frac{1+\sqrt{5}}{2}$ the golden ratio.
\end{coroll}
\begin{remark}
The result of the case $j=0$ above has already been obtained in~\cite{Schaeff1998,Left2000} using bijections with maps and 2-stack sortable permutations. We remark that the sequence of coefficients of $T_0(z)$ given by $1,1,2,6,22,91,408,\dots$ appear as sequence $A000139$ in the Online Encyclopedia of Integer Sequences~\cite{Sloane}. The coefficients of $T_1(z)$ given by $1,1,3,11,46,209,1006,\dots$ are not listed there.
\end{remark}

\begin{remark}
Banderier and Flajolet~\cite{BandFla2002} have studied directed lattice paths in the plane and pointed out the importance
of the so-called characteristic polynomial. For so-called simply paths, defined by the set of steps $\mathcal{S}=\{\mathbf{a}_1,\dots,\mathbf{a}_d\}$, with $\mathbf{a}_{\ell}=(1,b_{\ell})$ and $b_{\ell}\in\Z$, for $1\le \ell \le d$, the characteristic polynomial is given by $P(X)=\sum_{\ell=1}^d X^{b_{\ell}}$, and the characteristic equation of simple paths is given by $1-zP(X)=0$. 
Note that the equation for the series $X$ for embedded ternary trees with step set $\mathcal{S}=\{(1,-1),(1,0),(1,1)\}$ is given by $X=zT^2X(X^{-1}+1+X)$; for embedded binary trees~\cite{Bou2006} with step set $\mathcal{S}=\{(1,-1),(1,1)\}$ it is given by $X=zT(z)X(X^{-1} +X)$.
We conjecture that an equation of a similar type $1-z(T(z))^{|\mathcal{S}|-1}P(X)=0$, where $P(X)=\sum_{\ell=1}^{d}X^{b_\ell}$, is important at least for the analysis of embedded trees with a symmetric step set $\mathcal{S}$, where $T(z)=1+z(T(z))^{|S|}$ denotes the generating function of the number of embedded trees. 

\end{remark}

\begin{proof}
In order to proof the result above it is sufficient to check that
the series $T_j(z)$ satisfies Equation~\ref{BOUTT1} and the initial
condition $T_{-1}(z)=1$, which is a simple task. In order to
\emph{discover} the solution we use the method from Bouttier et
al.~\cite{Boutt2003}, see also Di Francesco~\cite{Francesco2005}. Following~\cite{Boutt2003} we use that fact
that for $j$ tending to infinity we have $T_j(z)\to T(z)$. Hence,
for $j$ tending to infinity one expects that $T_j(z)$ can be written
as $T_j(z)=T(1+o(1))$. We make the Ansatz $T_j(z)=T(1-\rho_j)$, with
$\rho_j\to 0$ as $j$ tends to infinity. We expend
Equation~\ref{BOUTT1} with respect to the Ansatz and compare the
terms tending at a similar rate to zero as $j$ tends to infinity. We
get
\begin{equation*}
\rho_j=zT^2(\rho_{j-1}+\rho_j+\rho_{j+1}).
\end{equation*}
We assume that $X=X(z)$ is a formal power series such that $|X(z)|<1$ for a given $r>0$ with $|z|<r< R$, where $R=4/27$ denotes the radius of convergence
of $T(z)$. An Ansatz $\rho_j=\alpha_1 X^j$ leads to the so-called characteristic equation
\begin{equation*}
X=zT^2(1+X+X^2),\quad\text{or equivalently}\quad 1=zT^2(X+1+\frac1X).
\end{equation*}
Note that the equation above is identical to the characteristic equation of Motzkin path $1=z(X+1+\frac1X)$ with 
$z$ being replaced by $zT^2$, where $\mathcal{S}=\{(1,1),(1,0),(1,-1)\}$, see~\cite{BandFla2002},. The equation above implies several other relations with respect to $zT^2=(T-1)/T$,
\begin{equation*}
\frac1X +X =\frac{1-zT^2}{zT^2}=\frac{1}{T-1}.
\end{equation*}
Moreover, the equations relating the series $T$ and $X$ of Theorem~\ref{DEEBtheBOUTT1} are easily derived in a similar manner.
Consequently, we define the series $X=X(z)$ as the solution of the equation above with $|X|<1$ and $X(0)=0$.
We make the improved Ansatz $\rho_j=\sum_{i\ge 1}\alpha_i(X^{j})^{i}$ and compare the terms with the same order of magnitude
in~\eqref{BOUTT1} as $j$ tends infinity.
Using the fact that $1/(zT^2)-1=X+1/X$ we obtain the system of recurrences
\begin{equation*}
\begin{split}
\alpha_{i+1}\Big(X^{i+1}+\frac1{X^{i+1}}-X-\frac1X\Big)&=\sum_{\ell=1}^{i}\alpha_{\ell}\alpha_{i+1-\ell}\big(\frac{1}{X^{\ell}}+X^{\ell}+X^{i+1-2\ell}\big)\\
&\quad -\sum_{\substack{\ell_1+\ell_2+\ell_3=i+1\\\ell_1,\ell_2,\ell_3\ge 1}}\alpha_{\ell_1}\alpha_{\ell_2}\alpha_{\ell_3}\frac{X^{\ell_3}}{X^{\ell_1}},\quad i\ge 0,
\end{split}
\end{equation*}
which determine $\alpha_{i+1}$ uniquely in terms of $\alpha_1$ and rational functions of $X$,
with $\alpha_1$ unspecified yet. Although it seems at first glance hopeless to obtain a closed form solution,
it is possible to guess a solution with the help of computer algebra software, i.e.~the help of \texttt{Maple},
which is easily verified using induction. We obtain the surprisingly simple solution
\begin{equation*}
\alpha_i=\frac{\alpha_1^i X^{i-1}(1-X^i)}{(1-X)(1-X)^{i-1}(1-X^2)^{i-1}},\quad i\ge 1.
\end{equation*}
By our initial assumption we have $|X|<1$ for $z$ small enough. We set
$\alpha_1=\alpha_1(X)=\lambda\cdot(1-X)(1-X^2)X^2$, where $\lambda\in \R$ being small enough and independent of $X$.
Consequently, we get after summation the solution
\begin{equation}
\label{DEEBgensolution}
T_j(z,\lambda)=T(z)\frac{(1-\lambda X^{j+2}(z))(1-\lambda X^{j+5}(z))}{(1-\lambda X^{j+3}(z))(1-\lambda X^{j+4}(z))},
\end{equation}
which satisfies the relation~\ref{BOUTT1} for all $j$, not necessarily larger than zero.
In order to see this, we use the notation $v_j=1-\lambda X^{j+1}$; the recurrence relation for $T_j(z)$, given in Lemma~\ref{BOUTT1}, 
implies that we have to show 
\begin{equation*}
T(z)v^2_{j+1}v_{j+2}v_{j+3}v^2_{j+4}=v_{j+1}v_{j+2}^2v_{j+3}^2v_{j+4}+zT^3(z)v_jv_{j+1}v_{j+2}v_{j+3}v_{j+4}v_{j+5},
\end{equation*}
which is easily seen to be true for all $j\in \Z$, with respect to $zT^3=T-1$ and $T=\frac{1+X+X^2}{1+X^2}$. Finally, adapting to the initial condition $T_{-1}(z)=1$ implies that $\lambda=1$.
\end{proof}

\subsection{The number of vertices with a given label}
Following~\cite{Bou2006} we are interested in the number of ternary trees of size $n$ with a given number
of vertices with label $j$. In order to treat this problem we introduce
a sequence of bivariate generating functions $S_j(z,u)$, where $z$ marks the number of internal nodes, and $u$ the number of nodes with label $j$,
with $j\in\Z$. Our first observations is already the key point, namely that due to the definition of the embedding of ternary trees, or in other words 
the symmetry of the step set $\mathcal{S}=\{(1,1),(1,0),(1,-1)\}$, we have the symmetry $S_j(z,u)=S_{-j}(z,u)$, for all $j\in\Z$; moreover we have $S_j(z,1)=T(z)$. The starting point of our considerations is the following lemma.
\begin{lemma}
\label{BOUTT2}
The series $S_j(z,u)$ satisfies
\begin{equation*}
\begin{split}
S_j(z,u)&=1+ zS_{j-1}(z,u)S_j(z,u)S_{j+1}(z,u),\quad \text{for}\,\, j\in\Z\setminus\{0\},\\
S_0(z,u)&=1+ uzS_{-1}(z,u)S_0(z,u)S_{1}(z,u)=1+uzS_0(z,u)S^2_{1}(z,u),\quad\text{for}\,\, j=0.
\end{split}
\end{equation*}
\end{lemma}

\begin{proof}
Using the arguments of~\cite{Bou2006} we observe that the series $S_j(z,u)$ is also the generating function of trees rooted at a node with label $j$, counted by the number of nodes labelled zero. The formal description sketched Figure~\ref{DEEBfig4} translates in the stated system of equations and the result follows.
\end{proof}

\begin{theorem}
\label{DEEBtheBOUTT2}
Let $S_j(z,u)$ be the generating function of ternary trees with counting the number of nodes with label $j$. Then $S_j(z,u)$ is given by the
the following expression
\begin{equation*}
S_j(z,u)=T(z)\frac{(1+\mu X^{j+1}(z))(1+\mu X^{j+4}(z))}{(1+\mu X^{j+2}(z))(1+\mu X^{j+3}(z))},\quad \text{for}\,\,j\ge -1,
\end{equation*}
where $T=T(z)$ is defined by~\eqref{DEEBteReqn0yi}, i.e.~$T=1+zT^3$, the series
$X=X(z)$ is specified in Theorem~\ref{DEEBtheBOUTT1}, and the power series $\mu=\mu(z,u)$ is defined as the unique formal power series with $\mu(z,1)=0$, satisfying the relations
\begin{equation*}
\mu=(u-1)\frac{(1+\mu X)(1+\mu X^2)^2(1+\mu X^5)}{(1+X)^2(1-X)^3(1-\mu^2X^5)}.
\end{equation*}
\end{theorem}

\begin{proof}
We have already seen in the proof of Theorem~\ref{DEEBtheBOUTT2},\eqref{DEEBgensolution} that the general solution of the set of equations for $j\neq 0 $ in Lemma~\ref{BOUTT2} is
given by
\begin{equation*}
T_j(z,\lambda)=T(z)\frac{(1-\lambda X^{j+2}(z))(1-\lambda X^{j+5}(z))}{(1-\lambda X^{j+3}(z))(1-\lambda X^{j+4}(z))}.
\end{equation*}
We have to determine $\lambda$ in such a way that the equation for $j=0$ in Lemma~\ref{BOUTT2} is satisfied.
Using the general solution stated above and the formulas expression $T(z)$ and $z$ as functions of $X=X(z)$ we obtain the equation
\begin{equation*}
\frac{(1+X+X^2)}{(1+X^2)}\frac{(1-\lambda X^{2}(z))(1-\lambda X^{5}(z))}{(1-\lambda X^{3}(z))(1-\lambda X^{4}(z))}
= 1+ u \frac{X}{(1+X^2)}\frac{(1-\lambda X^{2}(z))(1-\lambda X^{3}(z))(1-\lambda X^{6}(z))^2}{(1-\lambda X^{5}(z))(1-\lambda X^{4}(z))^3}.
\end{equation*}
Now we simply write $u=(u-1)+1$, and obtain after simple manipulations
\begin{equation*}
-\lambda X= (u-1)\frac{(1-\lambda X^2)(1-\lambda X^3)^2(1-\lambda X^6)}{(1+X)^2(1-X)^3(1-\lambda^2X^7)}
\end{equation*}
We set $\mu=-\lambda X$, and obtain the stated result.
\end{proof}

\subsection{The number of vertices with given labels\label{DEEBSubsectiongen1}}
In this section we derive a generalization of our previous result concerning the enumeration of embedded ternary trees of size $n$ with respect
to the number of nodes with label $j$. The crucial fact in the derivation of the previous result was the symmetry relation $S_j(z,u)=S_{-j}(z,u)$, which allowed to conclude the special equation for the case $j=0$. Hence, when generalizing the counting problem we have to take
care to preserve some symmetry in order to set up a system of suitable equations.
First, we are interested in counting two statistics at the same time,
namely the number of nodes with label $j$, and the number of nodes with label contained in $\{j-1,j+1\}:=\{j\pm 1\}$.
Let $S_j=S_{j,j\pm 1}(z,u_0,u_1)$ denote the generating function of ternary trees, where the variable $w_0$ counts the number of nodes
with label $j$ and $u_1$ counts the number of nodes with labels contained in $\{j\pm 1\}$. Since we both types of labels $j-1$ and $j+1$ are counted by the same variable, we have symmetry with respect to the vertical line $j$. More precisely, we obtain the following result.
\begin{lemma}
\label{BOUTT3}
The series $S_{j,j\pm 1}(z,u_0,u_1)$ satisfies the symmetry relation
$$S_{j,j\pm 1}(z,u_0,u_1)=S_{-j,-j\pm 1}(z,u_0,u_1)$$ for all $j\in\Z$.
Moreover, $S_j=S_{j,j\pm 1}(z,u_0,u_1)$ is determined by the system of equations
\begin{equation*}
\begin{split}
S_j&=1+ zS_{j-1}S_jS_{j+1},\quad \text{for}\,\, j \ge 2,\\
S_1&=1+ u_1zS_{0}S_{1}S_2,\quad\text{for}\,\, j=1,\\
S_0&=1+ u_0zS_{-1}S_0S_{1}=1+u_0 zS_0S^2_{1},\quad\text{for}\,\, j=0.
\end{split}
\end{equation*}
\end{lemma}
The proof is identical to the proof of Lemma~\ref{BOUTT2} and follows the argumentations of~\cite{Bou2006}. Therefore we omit it.
We obtain the following result.
\begin{theorem}
\label{DEEBtheBOUTT3}
The generating function $S_{j,j\pm 1}(z,u_0,u_1)$ of ternary trees, counting the number of nodes with label $j$ and the number of nodes with label in $\{j\pm 1\}$ is given by the
the following expression
\begin{equation*}
S_j(z,u)=T(z)\frac{(1+\mu X^{j+1}(z))(1+\mu X^{j+4}(z))}{(1+\mu X^{j+2}(z))(1+\mu X^{j+3}(z))},\quad \text{for}\,\,j\ge -1,
\end{equation*}
where $T=T(z)$ is defined by~\eqref{DEEBteReqn0yi}, i.e.~$T=1+zT^3$, the series
$X=X(z)$ is specified in Theorem~\ref{DEEBtheBOUTT1}, and the power series $\mu=\mu(z,u_0,u_1)$ is defined as the unique formal power series with $\mu(z,1,1)=0$, satisfying the relations
\begin{equation*}
\mu=(u_0-1)\frac{(1+\mu X)(1+\mu X^2)^2(1+\mu X^5)}{(1+X)^2(1-X)^3(1-\mu^2X^5)} + (u_1-1)\frac{(1+ X+X^2)(1+\mu X^2)^2(1+\mu X^3)^2}{X(1+X)^2(1-X)^3(1-\mu^2X^5)}.
\end{equation*}
\end{theorem}
\begin{remark}
Note that the series $\mu(z,u_0,1)$ reduces to the series stated in Theorem~\ref{DEEBtheBOUTT3}, as it should. Hence, Theorem~\ref{DEEBtheBOUTT3} is a generalization of our previous result. Furthermore, the series $S_{j,j\pm 1}(z,1,u_1)$, with $\mu(z,1,u_1)$, counts the number of nodes with label in $\{j \pm 1\}$.
\end{remark}
\begin{proof}
We proceed similar to the proof of Theorem~\ref{DEEBtheBOUTT2}. The equations for the series $S_1$ and $S_0$ given in Lemma~\ref{BOUTT3} imply the equations
\begin{equation*}
S_1=\frac{1}{1-w_1 zS_0S_2},\quad S_0=\frac{1}{1-w_0 zS_1^2}.
\end{equation*}
Consequently, we obtain by substituting the equation for $S_0$ into the equation for $S_1$ the result
\begin{equation*}
S_1= \cfrac{1}{1-\cfrac{w_1 zS_2}{1-w_0 zS_1^2}}.
\end{equation*}
As in the proof of Theorem~\ref{DEEBtheBOUTT2} we use the fact that the general solution of the set of equations for $j\neq \{-1,0,1\} $ in Lemma~\ref{BOUTT3} is given by
\begin{equation*}
T_j(z,\lambda)=T(z)\frac{(1-\lambda X^{j+2}(z))(1-\lambda X^{j+5}(z))}{(1-\lambda X^{j+3}(z))(1-\lambda X^{j+4}(z))}.
\end{equation*}
Simple manipulations leads to the stated result with respect to $\mu=-\lambda X$.
\end{proof}

Now discuss a more general counting problem. We are interested in the number of nodes with labels in $\{j\}=\{j \pm 0\}$, counted by the variable $u_0$, labels in $\{j \pm 1\}$, counted by $u_1$, up to labels in $\{j \pm m\}$, counted by $u_m$. Note that the cases $m=0$ and $m=1$ corresponds to the counting problems treated in Theorems~\ref{DEEBtheBOUTT2},\ref{DEEBtheBOUTT3}. We use the vector notation $\mathbf{u}=(u_0,\dots,u_m)$.
Let $S_j(z,\mathbf{u})=S_{j\pm 0,j\pm 1,\dots,j\pm m}(z,u_0,\dots,u_m)$ denote the generating function of ternary trees, where the variable $u_{\ell}$ counts the number of nodes with label in $\{j \pm l\}$, for $0\le \ell \le m$.
\begin{lemma}
\label{BOUTT4}
The series $S_j=S_j(z,\mathbf{u})$ satisfies the symmetry relation
$S_{j}=S_{-j}$ for all $j\in\Z$.
Moreover, $S_j$ is determined by the system of equations
\begin{equation*}
\begin{split}
S_j&=1+ zS_{j-1}S_jS_{j+1},\quad \text{for}\,\, j \ge m+1,\\
S_{j}&=1+ u_{j}zS_{j-1}S_{j}S_{j+1},\quad\text{for}\,\, 0\le j\le m,\\
S_0&=1+ u_0zS_{-1}S_0S_{1}=1+u_0 zS_0S^2_{1},\quad\text{for}\,\, j=0.
\end{split}
\end{equation*}
\end{lemma}
The equations above are obtained in the same manner as the equations in the Lemmata~\ref{BOUTT2},\,\ref{BOUTT3}. Let $<a_n,\dots,a_1,a_0>_n$ denote the finite continued fraction in the variables $a_n,\dots,a_0$,
\begin{equation}
\label{DEEBcfrac}
<a_n,\dots,a_1,a_0>_n=\cfrac{1}{1-\cfrac{a_n}{\ddots-\cfrac{\vdots}{1-a_0}}}
\end{equation}
Consequently, the system of equations stated in Lemma~\ref{BOUTT4}
for $0\le j\le m$ implies, with respect to the notation above, the
equation
\begin{equation}
\label{DEEBgen1}
S_m=<zu_mS_{m+1},zu_{m-1}S_{m}\dots,zu_1S_{2},zu_0 S_{1}^2>_m.
\end{equation}

We will use the following lemma.
\begin{lemma}
\label{BOUTT5}
The finite continued fraction $<a_n,\dots,a_1,a_0>_n$, defined in~\eqref{DEEBcfrac}, can be written in terms of the family of polynomials
$k_n(a_n,\dots,a_0)$ in the following way,
\begin{equation*}
<a_n,\dots,a_1,a_0>_n=\frac{k_{n-1}(a_{n-1},\dots,a_0)}{k_n(a_n,\dots,a_0)}\quad\text{with}\,\, k_n(a_n,\dots,a_0)=1+\sum_{\ell=1}^{\lfloor\frac{n}{2}\rfloor+1}(-1)^{\ell}\sum_{\substack{(i_1,\dots,i_{\ell})\\i_{j-1}+2\le i_j\le n\\0\le j\le \ell, \,i_{-1}:=-2}}a_{i_1}\dots a_{i_{\ell}}.
\end{equation*}
\end{lemma}
\begin{proof}
The finite continued fraction $<a_n,\dots,a_1,a_0>_n$ can be written as a quotient of the polynomials $h_n=h_n(a_n,\dots,a_0)$ and $k_n=k_n(a_n,\dots,a_0)$,
\begin{equation*}
<a_n,\dots,a_1,a_0>_n=\frac{h_n}{k_n},\quad\text{with}\,\, h_0=1,\,k_0=1-a_0,\quad\text{and we define}\,\,k_{-1}=1.
\end{equation*}
We easily obtain the recurrence relations
\begin{equation*}
h_n=k_{n-1},\quad\text{and}\,\, k_n=k_{n-1}-a_n k_{n-2},\quad\text{for}\,\, n\ge 1,
\end{equation*}
with respect to the initial values of $k_n$ stated above.
Now we easily prove with induction that the polynomials $k_n$ are given by the stated explicit expression.
\end{proof}
Using again the general solution of the set of equations for $j\neq \{-m,\dots,m\}$ given by~\eqref{DEEBgensolution}
we obtain from~\eqref{DEEBgen1} an equation for the unknown series in terms of $\lambda=\lambda(z,\mathbf{u})$,
\begin{equation*}
\begin{split}
&S_m\cdot k_m(zu_mS_{m+1},zu_{m-1}S_{m}\dots,zu_1S_{2},zu_0 S_{1}^2)-k_{m-1}(zu_{m-1}S_{m}\dots,zu_1S_{2},zu_0 S_{1}^2)=0,\\
&\text{with}\,\,S_j=S_j(z,\mathbf{u})=T(z)\frac{(1-\lambda X^{j+2}(z))(1-\lambda X^{j+5}(z))}{(1-\lambda X^{j+3}(z))(1-\lambda X^{j+4}(z))},
\end{split}
\end{equation*}
where the polynomials $k_n=k_n(a_n,\dots,a_0)$ are given in Lemma~\ref{BOUTT5}.

One can easily check that for small $m=1,2,3,4$ one curiously always obtains a ``simple'' equation for $\mu=\mu(z,\mathbf{u})=-\lambda X$ of the form
\begin{equation*}
\mu=\sum_{j=0}^{m}(u_j-1)\cdot f_{j,m}(\mu,X)+\delta_m(\mathbf{u}\boldsymbol{-1},\mu,X),
\end{equation*}
where $\delta_m$ is a polynomial in the $u_j-1$, $0\le j\le m$ with degree greater than $1$. Unfortunately, we have not been able to prove this closed form expression for general $m$.

\begin{remark}
Let $S_j(z,\mathbf{u})=S_{j\pm 0,j\pm 1,\dots,j\pm m}(z,u_0,\dots,u_m)$ denote the generating function of \emph{embedded binary trees}, as treated in~\cite{Bou2006}, where the variable $u_{\ell}$ counts the number of nodes with label in $\{j \pm l\}$, for $0\le \ell \le m$.
Using the considerations of~\cite{Bou2006} and proceeding similar to the arguments of the ternary tree case, one may obtain a similar description for the generating function $S_j(z,\mathbf{u})$ in
terms of the formal power series $\lambda=\lambda(z,\mathbf{u})$,
\begin{equation*}
S_j(z,\mathbf{u})=T(z)\frac{(1-\lambda X^{j+2}(z))(1-\lambda X^{j+7}(z))}{(1-\lambda X^{j+4}(z))(1-\lambda X^{j+4}(5))}, \quad T(z)=1+zT^2(z)=\frac{1+X^2}{1-X+X^2},\quad ,
\end{equation*}
with the series $T=T(z)$ denoting the generating of binary trees $T=1+zT^2$, the series $X$ defined by
$X=zT X(X^{-1}+X)$ satisfying $X(0)=0$, and the series $\lambda=\lambda(z,\mathbf{u})$ being defined by
\begin{equation*}
S_m - \sum_{i=1}^{m+1}\prod_{\ell=i}^{n}(zu_{\ell}S_{\ell+1}) - zu_0S_1^2\prod_{\ell=1}^{n}(zu_{\ell}S_{\ell+1})=0.
\end{equation*}
The case $m=0$ simplifies to a result stated in~\cite{Bou2006}. As for ternary trees, we expect
a ``simple'' equation for $\mu=\mu(z,\mathbf{u})=-\lambda X^2$ of the form
\begin{equation*}
\mu=\sum_{j=0}^{m}(u_j-1)\cdot f_{j,m}(\mu,X)+\delta_m(\mathbf{u}\boldsymbol{-1},\mu,X),
\end{equation*}
where $\delta_m$ is a polynomial in the $u_j-1$, $0\le j\le m$ with degree greater than $1$.
Again, this has been checked for small values of $m$, but still lacks a proof for general $m$.
\end{remark}

\begin{remark}
It should be possible to refine the stated results to refine the stated results 
by seperating the contributions of nodes $-j$ and $j$, $1\le j \le m$, i.e.~by introducing the refined generating function
$S_j(z,\mathbf{u})=S_{j\pm 0,j\pm 1,\dots,j\pm m}(z,u_0,u_1,u_{-1},\dots,u_m,u_{-m})$ and exploiting the symmetry
$S_{j\pm 0,j\pm 1,\dots,j\pm m}(z,u_0,u_1,u_{-1},\dots,u_m,u_{-m})=S_{-j\pm 0,-j\pm 1,\dots,-j\pm m}(z,u_0,u_{-1},u_1\dots,u_{-m},u_m)$.
\end{remark}

\section{The distribution of depths of a leaf in ternary trees\label{DEEBTERleaves}}
We obtain the following result for the distribution of types of depths of the external node (leaf) $s$, where
the leaves of the ternary tree are enumerated from the left to the right.
\begin{theorem}
\label{DEEBtheAlois}
The number $T_{n,s,m_1,m_2,m_3}$ of ternary increasing trees of size $n$ where leaf $s$ has $m_i$ edges of type $e_i$, $i=1,2,3$, is given for $n\ge 1$ by the following explicit formulas.
\begin{equation*}
T_{n,s,m_1,m_2,m_3}=
\begin{cases}
\displaystyle{\frac{(2m_3+m_2)(2m_1+m_2)\binom{m_1+m_2+m_3}{m_1,m_2,m_3}\binom{3s_1+s_2-m_3-\mu_2}{s_1-m_3-\mu_2}\binom{3n-m_1-2m_2-3s_1+3\mu_2}{n-m_1-m_2-s_1+\mu_2}}{(3s_1+s_2-m_3-\mu_2)(3n-m_1-2m_2-3s_1+3\mu_2)}},\\
\quad\text{for}\quad s=2s_1+s_2,\,\, 0\le s_1 \le n\,\,\text{and}\,\,s_2\in\{0,1\},\quad\\
\quad\text{with}\quad 0\le s_1\le m_3, m_2=s_2+2\mu_2,\,\,0\le \mu_2 \le s_1-m_3,\\
\quad\text{and}\quad(m_1,m_2),(m_3,m_2)\neq (0,0);\\[0.2cm]
\displaystyle{\frac{2m_1}{3n-m_1}\binom{3n-m_1}{2n}},\quad\text{for}\,\,(m_2,m_3)=(0,0),\quad s=0,\quad 0\le m_1\le n,\\[0.4cm]
\displaystyle{\frac{2m_3}{3n-m_3}\binom{3n-m_3}{2n}},\quad\text{for}\,\,(m_1,m_2)=(0,0),\quad s=2n,\quad 0\le m_3\le n,\\[0.2cm]
0\qquad\text{for all other cases}.
\end{cases}
\end{equation*}
\end{theorem}
\begin{remark}
Under the random tree model, i.e.~every ternary tree of size $n$ is equally likely, we can consider the
random variable $H_{n,s}^{[i]}$ counting the number of edges of type $i$, with $1\le i\le 3$, on the path from the external node numbered $s$
to the root. The joint distribution can be expression in terms of the numbers $T_{n,s,m_1,m_2,m_3}$ and the total number of ternary trees of size $n$,
\begin{equation}
\P\{H_{n,s}^{[1]}=m_1,H_{n,j}^{[2]}=m_2,H_{n,s}^{[3]}=m_3\}=\frac{T_{n,s,m_1,m_2,m_3}}{T_{n}}
\end{equation}
Hence, one may derive limiting distribution results in the manner of~\cite{Alois2009} from the explicit expression for $T_{n,s,m_1,m_2,m_3}$ using Stirling's formula.
\end{remark}
\begin{proof}
Our methods are based on the considerations of Panholzer~\cite{Alois2009}, see also the works of Panholzer and Prodinger~\cite{PP1997,PP2002}.
Let $\mathbf{v}=(v_1,v_2,v_3)$ where the variable $v_i$  counts the number of edges of type $i$, with $i=1,2,3$, or equivalently the $v_1$ denotes a left edge $e_1$,
$v_2$ a center edge $e_2$, and $v_3$ a right edge $e_3$. Subsequently we will use the multiindex notation $\mathbf{m}=(m_1,m_2,m_3)$ and $\mathbf{v}^{\mathbf{m}}=v_1^{m_1}v_2^{m_2}v_3^{m_3}$; moreover
$\mathbf{m}\ge 0$ should be interpreted componentwise. We introduce the generating function of $T_{n,s,\mathbf{m}}$, the number of ternary trees
of size $n$ where leaf $j$ has $m_1$ left steps, $m_2$ center steps, and $m_3$ right steps, with $m_1+m_2+m_3\le n$, with $m_i\ge 0$.
\begin{equation*}
F(z,u,\mathbf{v})=\sum_{n\ge 0}\sum_{s=0}^{2n}\sum_{\mathbf{m}\ge 0}T_{n,s,\mathbf{m}}z^nu^s\mathbf{v}^{\mathbf{m}}.
\end{equation*}
Note that since we enumerate the leaves from left to right, each
center step increases the label of the leaf by one and each right
step increases the label of the leaf by two. Moreover each new
internal node increases the number of leaves by two. Hence, we get
the additional condition that for the $s$-th leaf $s-2m_3-m_2\equiv
0\mod 2$. The formal description of ternary trees translates into
the following functional equation for $F(z,u,\mathbf{v})$, where
$T(z)$ denotes the generating function of ternary trees.
\begin{equation*}
F(z,u,\mathbf{v})=1 + zv_1F(z,u,\mathbf{v}) + zuv_2T(zu^2)F(z,u,\mathbf{v})T(z)+zu^2v_3\big(T(zu^2)\big)^2F(z,u,\mathbf{v}).
\end{equation*}
We obtain
\begin{equation*}
\begin{split}
F(z,u,\mathbf{v})&= \frac{1}{1-zv_1-zuv_2T(zu^2)T(z)-zu^2v_3\big(T(zu^2)\big)^2}\\
&=\frac{1}{(1-zv_1)\Big(1-\cfrac{zuv_2T(zu^2)T(z)}{1-zv_1}\Big)\Big(1-\cfrac{v_3zu^2\big(T(zu^2)\big)^2}{\Big(1-\cfrac{zuv_2T(zu^2)T(z)}{1-zv_1}\Big)(1-z_1v_1)}\Big)},
\end{split}
\end{equation*}
and consequently
\begin{equation*}
\begin{split}
[z^nu^s\mathbf{v}^{\mathbf{m}}]F(z,u,\mathbf{v})&=[z^nu^s]\binom{m_1+m_2+m_3}{m_1,m_2,m_3}z^{m_1+m_2+m_3}u^{2m_3+m_2}\big(T(zu^2)\big)^{2m_3+m_2}
\big(T(z)\big)^{2m_1+m_2}\\
&=[z^{n-m_1-m_2-m_3}u^{s-2m_3-m_2}]\binom{m_1+m_2+m_3}{m_1,m_2,m_3}\big(T(zu^2)\big)^{2m_3+m_2}
\big(T(z)\big)^{2m_1+m_2}.
\end{split}
\end{equation*}
Extraction of coefficients will be carried out by an application of our previous result stated in Equation~\ref{DEEBteReqn4};
We assume that the number of the leaf $s$, with $0\le s\le 2n$, is given as $s=2s_1+s_2$, for $0\le s_1 \le n$ and $s_2\in\{0,1\}$,
with $s_2=0$ for $s_1=n$. Consequently, we can rewrite the condition $s-2m_3-m_2\equiv 0\mod 2$ into
$0\le m_3\le s_1$ and $m_2=s_2+2\mu_2$, with $0\le \mu_2 \le s_1-m_3$. First we assume that both $(m_1,m_2)\neq (0,0)$ and $(m_2,m_3)\neq (0,0)$.
Using~\eqref{DEEBteReqn4} we get
\begin{equation*}
\begin{split}
&[z^{n-m_1-m_2-m_3}u^{s-2m_3-m_2}]\big(T(zu^2)\big)^{2m_3+m_2}\big(T(z)\big)^{2m_1+m_2}=\\
&\qquad\quad=\frac{(2m_3+m_2)(2m_1+m_2)\binom{3s_1+s_2-m_3-\mu_2}{s_1-m_3-\mu_2}\binom{3n-m_1-2m_2-3s_1+3\mu_2}{n-m_1-m_2-js_1+\mu_2}}{(3s_1+s_2-m_3-\mu_2)(3n-m_1-2m_2-3s_1+3\mu_2)}.
\end{split}
\end{equation*}
Note that we can write the expression above solely in $n,s,m_1,m_2,m_3$,
\begin{equation*}
\begin{split}
&[z^{n-m_1-m_2-m_3}u^{s-2m_3-m_2}]\big(T(zu^2)\big)^{2m_3+m_2}\big(T(z)\big)^{2m_1+m_2}=\\
&=\qquad\quad\frac{(2m_3+m_2)(2m_1+m_2)\binom{s-m_3+\frac{s-m_2}{2}}{\frac{s-m_2}{2}-m_3}\binom{n-s-m_1-\frac{s+m_2}{2}}{n-m_1-\frac{s+m_2}{2}}}{(s-m_3+\frac{s-m_2}{2})(n-s-m_1-\frac{s+m_2}{2})}.
\end{split}
\end{equation*}
Second, we assume that $(m_2,m_3)=(0,0)$, which implies that $s=0$ and $0\le m_1\le n$. We obtain
\begin{equation*}
[z^nu^j\mathbf{v}^{\mathbf{m}}]F(z,u,\mathbf{v})=[z^{n-m_1}]\big(T(z)\big)^{2m_1}=\frac{2m_1}{3n-m_1}\binom{3n-m_1}{2n}.
\end{equation*}
Third, lets assume that $(m_1,m_2)=(0,0)$, which implies that $s=2n$ and $0\le m_3\le n$
We obtain
\begin{equation*}
[z^nu^j\mathbf{v}^{\mathbf{m}}]F(z,u,\mathbf{v})=[z^{n-m_3}u^{s-2m_3}]\big(T(zu^2)\big)^{2m_3}=\frac{2m_3}{3n-m_3}\binom{3n-m_3}{2n},
\end{equation*}
which finishes the proof of Theorem~\ref{DEEBtheAlois}.
\end{proof}

\section{Outlook: Embedded d-ary trees}
We discuss generalizations of the results obtained in this note for embedded ternary trees to embedded $d$-ary trees, $d\ge 2$, with respect to the natural embedding of $d$-ary trees defined in Subsection~\ref{DEEBsecembedd}.

\subsection{The distribution of depths of a leaf in d-ary trees}
Subsequently, we will show that the study of the parameter depths of external nodes (leaf) $s$, with $0\le s\le dn$, where the $d$ different types
of depths stem from the embedding and the number of the leafs are due to an inorder-traversal, can be treated in the same manner as in the ternary (and binary) case.
We introduce the generating function
\begin{equation*}
F(z,u,\mathbf{v})=\sum_{n\ge 0}\sum_{s=0}^{dn}\sum_{\mathbf{m}\ge 0}T_{n,s,\mathbf{m}}z^nu^s\mathbf{v}^{\mathbf{m}}.
\end{equation*}
As before we enumerate the leaves from left to right, each step (or equivalently an edge) of type $\ell$ increases the label of the leaf by $\ell-1$,
for $1\le \ell \le d$. Moreover, since each new internal node increases the number of leaves by $d-1$, we obtain the additional condition that for the $s$-th leaf $s-\sum_{\ell=1}^{d}m_{\ell}(l-1)\equiv 0\mod (d-1)$. The formal description of $d$-ary trees tranlates into the following functional equation for $F(z,u,\mathbf{v})$, where $T(z)$ denotes the generating function of $d$-ary trees.
\begin{equation*}
F(z,u,\mathbf{v})=1 + z\sum_{\ell=1}^{d}v_{\ell}u^{\ell-1}\big(T(zu^{d-1})\big)^{\ell-1}F(z,u,\mathbf{v})\big(T(z)\big)^{d-l}.
\end{equation*}
Consequently, we obtain
\begin{equation}
\label{DEEBd1}
F(z,u,\mathbf{v})=\frac{1}{1-z\sum_{\ell=1}^{d}v_{\ell}u^{\ell-1}\big(T(zu^2)\big)^{\ell-1}\big(T(z)\big)^{d-l}}.
\end{equation}
In order to extract coefficients we use the the following folklore result.
\begin{lemma}
\label{LEMext1}
The formal power series $1/(1-\sum_{\ell=1}^{k}x_{\ell}\alpha_{\ell})$ in the $k$ variables $x_1,\dots,x_k$
satisfies the expansion
\begin{equation*}
[x_k^{m_k}x_{k-1}^{m_{k-1}}\dots x_1^{m_1}]\frac{1}{1-\sum_{\ell=1}^{k}x_k\alpha_k}=
\binom{m_1+\dots+m_k}{m_1,\dots,m_k}\prod_{\ell=1}^{k}\alpha_{\ell}^{m_{\ell}}.
\end{equation*}
\end{lemma}
\begin{proof}
Observe that by the multinomial theorem the generating function of $\binom{m_1+\dots+m_k}{m_1,\dots,m_k}\prod_{\ell=1}^{k}\alpha_{\ell}^{m_{\ell}}$
is given by $1/(1-\sum_{\ell=1}^{k}x_{\ell}\alpha_{\ell})$, or use
\begin{equation*}
[x_k^{m_k}]\frac{1}{(1-\sum_{\ell=1}^{k}x_{\ell}\alpha_{\ell})^j}=
[x_k^{m_k}]\frac{1}{(1-\sum_{\ell=1}^{k-1}x_{\ell}\alpha_{\ell})^j\big(1-\frac{x_k\alpha_k}{(1-\sum_{\ell=1}^{k-1}x_{\ell}\alpha_{\ell})}\big)^j}=
\alpha_k^{m_k}\frac{\binom{m_{k}+j-1}{m_k}}{(1-\sum_{\ell=1}^{k}x_{\ell}\alpha_{\ell})^{m_k+j}}
\end{equation*}
and induction with respect to $k$ concerning the other variables $x_{k-1},\dots,x_1$.
\end{proof}
An application of Lemma~\ref{LEMext1} to~\eqref{DEEBd1} immediately provides
\begin{equation*}
\begin{split}
[z^nu^s\mathbf{v}^{\mathbf{m}}]F(z,u,\mathbf{v})&=\binom{m_1+\dots+m_d}{m_1,\dots,m_d}[z^nu^s]\prod_{\ell=1}^{d}\Big(zu^{\ell-1}\big(T(zu^{d-1})\big)^{\ell-1}\big(T(z)\big)^{d-\ell}\Big)^{m_l}\\
&=\binom{m_1+\dots+m_d}{m_1,\dots,m_d}[z^nu^s]z^{M_1}u^{M_2-M_1}\big(T(zu^{d-1})\big)^{M_2-M_1}\big(T(z)\big)^{dM_1-M_2},
\end{split}
\end{equation*}
with respect to the notations $M_1:=\sum_{\ell=1}^{d}m_{\ell}$ and $M_2:=\sum_{\ell=1}^{d}\ell m_{\ell}$.
Extraction of coefficients can be carried out via the correspondence
to the family of simply generated $d$-ary trees, i.e.~ the formula
\begin{equation*}
[z^n]\big(T(z)\big)^{k}=\frac{k}{dn+k}\binom{dn+k}{n},\quad\text{for}\,\, n\ge 0,\,\,k\ge 1,
\quad\text{with}\,\, T(z)=1+z\big(T(z)\big)^{d},
\end{equation*}
where $T(z)$ denotes the generating function of $d$-ary trees.

\subsection{Nodes with small labels}
Concerning the enumeration of trees with small labels, the case of $d\ge 4$ is much more difficult than the cases $d=2,3$. In the general case $d\ge 4$ the approach of Bouttier et al.~\cite{Boutt2003} is still applicable; one obtains the following system of equations for the generating function $T_j(z)$ of 
$(2d+1)$-ary tree with no label greater than $j$, 
\begin{equation*}
T_j(z)=1+z\prod_{\ell=-d}^{d}T_{\ell}(z),\quad j\ge 0,\quad\text{with initial values}\,\,T_j(z)=1,\quad\text{for}\,\,-d\le j \le 1;
\end{equation*}
the recurrence relations for $(2d)$-ary trees are identical with respect to skipping the index $k=0$ in the product above.
The Ansatz $T_j(z)=T(z)(1-\rho_j)$, with $\rho_j=\alpha \cdot X^j$, with an a priori unknown series $X=X(z)$, leads again to a characteristic equations of the type
\begin{equation*}
1-zT^{2d}(z)\sum_{\ell=-d}^{d}X^{\ell}=0,\quad\text{and}\quad 1-zT^{2d-1}(z)\sum_{\substack{-d\le \ell \le d\\l\neq 0}}X^{\ell}=0.
\end{equation*}
The resulting characteristic equation has more than one solution, see~\cite{BandFla2002} for a related problem, and consequently the refined Ansatz gets more involved. However, it seems possible to extend the results obtained for binary trees~\cite{Bou2006} and ternary trees at least to quaternary and quinary trees, as indicated by the results of~\cite{Boutt2003}. For example, one can obtain a
family of formal solutions for embedded $(2d+1)$-ary trees
\begin{equation*}
T_j(z)=T(z)\frac{(1-\lambda X^{d+1+j}(z))(1-\lambda X^{2d+3+j}(z))}{(1-\lambda X^{d+2+j}(z))(1-\lambda X^{2d+2+j}(z))},
\end{equation*}
where $X$ is a solution of the characteristic equation. Unfortunately, none of the roots of the characteristic equation
is sufficient for solving the counting problem; we believe that the proper solution $T_j(z)$ can be expressed in terms of all ``small roots'' of the characteristic equation stated above (in the sense of~\cite{BandFla2002}). The author is currently investigating into this matter.

Furthermore, concerning the general setting of embedded tree, see Subsection~\ref{DEEBform}, specified by the set of vectors $\mathcal{S}=\{\mathbf{a}_1,\dots,\mathbf{a}_k\}$, with $\mathbf{a}_{\ell}=(1,b_{\ell})$, we obtain the characteristic equation
\begin{equation*}
1-zT^{|\mathcal{S}|-1}(z)P(X)=0,\quad\text{with}\,\,P(X)=\sum_{\ell=1}^{k}X^{b_\ell}.
\end{equation*}
\subsection{Some open problems}
We did not manage to prove the conjectured form of the formulas for the generating function counting several labels at once. 
Moreover, it would be nice to extend the results for the generating functions of embedded binary trees~\cite{Bou2006} and ternary trees 
to two constraints, i.e.~asking for the number of trees with labels only in say $\{0,1,\dots, m\}$, $m\ge 0$. We refer the interested reader to the work of Bouttier et al.~\cite{Boutt2003II}, where a similar problem has been solved for families of plane trees using a $q$-theta function Ansatz, revealing an interesting connection to elliptic functions.

\section*{Acknowledgements}
The author wants to thank Alois Panholzer for valuable and fruitful discussions. Moreover, he thanks 
Mireille Bousquet-M\'elou for encouraging remarks and pointing out several errors.

\end{document}